\documentclass[12pt]{article}
\usepackage{amsmath,amssymb}
\usepackage{graphicx}
\usepackage{amsthm}

\def\RC{\mathcal{R}}

\def\R{\mathbf{R}}

\def\1{\mathbf{1}}

\def\al{\alpha}

\newtheorem{theorem}{Theorem}[section]

\newtheorem{remark}{Remark}

\newcommand{\la}{\lambda}

\newcommand{\om}{\omega}

\newcommand{\La}{\Lambda}

\begin{document}
\title{On the control over the distribution of ticks based on the extensions of the KISS model}

\author{Vassili N. Kolokoltsov\thanks{Department of Statistics, University of Warwick,
 Coventry CV4 7AL UK, and Higher School of Economics Moscow,
 ORCID 0000-0001-6574-4615,
  Email: v.kolokoltsov@warwick.ac.uk}}
\maketitle

\begin{abstract}
Ticks and tick-borne diseases present a well known threat to the health of people in many parts of the globe.
The scientific literature devoted both to field observations and to modeling the propagation of
ticks continues to grow. So far the majority of the mathematical studies were devoted to models
based on ordinary differential equations, where spatial variability was taken into account by a discrete parameter.
Only few papers use spatially nontrivial diffusion models, and they are devoted mostly to spatially homogeneous
equilibria. Here we develop diffusion models for the propagation of ticks stressing spatial heterogeneity. This
allows us to assess the sizes of control zones that can be created (using various available techniques) to produce
a patchy territory, on which ticks will be eventually eradicated. Using averaged parameters taken from various
field observations we apply our theoretical results to the concrete cases of the lone star ticks of North America
and of the taiga ticks of Russia.
\end{abstract}

\vspace{3mm}


\section{Introduction}

Ticks and tick-borne diseases present a well known threat to the health of people in many parts of the globe.
The scientific literature devoted both to field observations and to modeling the propagation of
ticks continues to grow.

So far the majority of the mathematical studies were devoted to models
based on ordinary differential equations, where spatial variability was taken into account by a discrete parameter.
One of the most basic models was developed in \cite{Gaff07}, which is an ODE in 4 parameters (number of ticks $V$, number of hosts $N$,
number of infected ticks $X$, number of infected hosts $Y$), and where the standard quadratic term of logistic models (expressing the internal competition of species) was replaced by a term limiting the growth of ticks by the number of available hosts.
Apart from the analysis of equilibria (and their stability) performed for this model, paper \cite{Gaff07} also suggests
a variant with multiple patches (with some fixed rates of transmission between them), however with only numerical analysis
for this advanced variant.

In \cite{White19} this model was further developed to include the dynamics of two pathogens.
Very complicated extensions (with various variables and parameters) were suggested and analysed in
\cite{Swit16}, \cite{Zhang21} and \cite{Nah21}. In \cite{Carv15} a fractional version of tick-host dynamics was developed.

A stochastic version (Markov-chain model) of the model from \cite{Gaff07} was developed in \cite{Mal17},
and the extinction predictions for deterministic and stochastic models were compared.
In \cite{Gaff11} the models from \cite{Gaff07} were enhanced by introducing some control
parameters and the related optimal control policies (for quadratic costs) were derived
aimed at reducing the propagation of diseases.

In \cite{Kash20} and \cite{Shu20} models with delays were suggested and analysed via various methods
(Lyapunov functions, Hopf bifurcations). A cellular-automaton model for the spatial propagation of
ticks was developed and analysed numerically in \cite{Vshiv13}.

In \cite{Mwa00} a stage structured population model was proposed and analysed. In the present
paper we extend partially this model by including spatial dynamics.

Paper \cite{Tosato21} investigates how different types of host targeted
tick-control strategies (combined acaricide-repellent control) affect tick population and disease
transmission, leading to some remarkable counterintuitive conclusions.

A reaction- diffusion model for the propagation of lyme disease was suggested in \cite{Caraco02}
extending a similar spatially trivial model of \cite{Caraco98}.
In \cite{Caraco02} the dynamics involved 9 density-variables: susceptible mice $M$, infected mice $m$,
questing larvae $L$, larvae infesting susceptible mice $V$, larvae infesting pathogen infected mice $v$,
susceptible questing nymphs $N$, questing infectious nymphs $n$,
uninfected adult ticks $A$, pathogen infected adult ticks $a$. Position-independent equilibria were obtained
and their local stability investigated.
This model was further developed in \cite{Zhao12}, where the global stability of spatially trivial equilibria was analysed
for the case of bounded domains and the travel wave solutions were constructed for the case of unbounded domains.

The exploitation of partial differential equations for the analysis of spatial distribution of various species
is a standard tool in ecology, see e.g.  \cite{Holmes94} and \cite{Holmes93}. The latter paper, in particular,
discusses the advantages and disadvantages of using the standard diffusion models, as compared with the models
based on the telegraph equation. One of the standard tool is the so-called KISS model (from Kierstead, Slobodkin
 and Skellam, see \cite{KiS53} and \cite{Skellam51}).
It was developed initially in the framework of the plankton propagation. It describes the minimal size of a patch,
where the population can survive under the killing (mathematically Dirichlet) conditions on the boundary.
For instance, in one dimensional case we look at the equation $\dot y=ay''+\la y $, with $a,\la>0$
(see \eqref{eqtickdiff1} below for all notations), on an interval $[0,R]$ with the
conditions $y(0)=y(R)=0$. The minimal $R$, for which a solution will not necessarily  tend to zero, as time goes
to infinity, is referred to as the critical patch size or the KISS scale. More generally one considers the behavior
of the population when the interval $[0,R]$ is surrounded by an infinite  non-beneficial territory, where the dynamics is described
by the equation $\dot y=ay''-\mu y$, with $\mu>0$. The critical size then decreases depending on the ratio $\la/\mu$. A review
of these results and their various extensions can be found in \cite{Okubo84}, see also \cite{Mu03}.

Let us mention a very important development of the KISS model for the description of several interacting species
(though less relevant to the present paper), in particular diffusive versions of the Lotka-Volterrra equations,
which dealt with patchiness (spatially nontrivial stationary solutions) and its stability,
see e.g.  \cite{Baron20}, \cite{Gourley96}, \cite{Mimura78},

The KISS  model is the starting point for our analysis. But eventually we are interested in creating control
zones (kind of barriers) with non-beneficial conditions that can stop the propagation of ticks. We first
extend the KISS model to the case when the non-beneficial region is bounded, interpreting its size as an additional
control parameter, and the diffusion coefficient is discontinuous at the interface of beneficial
and non-beneficial regions together with the mortality coefficient. Then we consider the modifications of the model
with various boundary conditions, most importantly with the periodic boundary conditions.
 Being seemingly not much relevant from the practical point of view, they provide
the key intermediate step to tackle our main questions: how to organise the sizes and properties of the regular patches of
(intermittent) control zones that would prevent the growth of the population of ticks.
Our main results are partial answers to this question.
We give first rather precise results for a homogeneous population, and then extend them partially
 to the structured populations, with several stages,
like larvae, nymphs and adults in ticks.

We suppose that in the control zone one is able to achieve the required level of the death rates or to
decrease the reproduction rates of ticks. In literature one can find various methods that can be used locally
for this purpose. It can be achieved by the treatment of the soil (planned burns, cleaning and
decreasing humidity levels, chemical treatment, see e.g. \cite{Tosato21}, \cite{Filippova85} and references therein),
 by tick-targeted strategies such as TickBots (see e.g. \cite{Gaff15} and \cite{Tosato21}), by specific actions on hosts,
like dipping of cattle with acaricide (see e.g. \cite{Ostfeld}  and \cite{Mount93}) or by biological methods
like infesting woods with ants (see \cite{HarrisAnts} for the related experiments on ants in USA and
\cite{Vshiv09} for the experiments in Russia) or the natural enemies of ticks like the insects Ixodiphagus hookeri
(see e.g. \cite{Filippova85}).

From the mathematical point of view, what we are doing here is getting criteria for negativity of certain differential operators,
or more precisely, looking for the precise estimates to their highest eigenvalue. This is one of the basic problems in the theory
of pseudo-differential equations. Many general results are available, which are often given, however, for the operators with smooth
(at least major) symbol and/or up to certain constants, see e.g. \cite{Feffer} and \cite{Sweezy93}. Here we solve this problem
 for a class of second order operators with piecewise continuous symbols. The negativity of the diffusion operators that govern
 the propagation of ticks ensures that the population will eventually die out (in the region considered).

The content of the paper is as follows.
In Section \ref{sectheory} our models are introduced and the main mathematical results are formulated.
The first introductory subsection deals with a simple scalar model of one-dimensional diffusion. We give the necessary
and sufficient conditions for the negativity of the corresponding diffusion operators
(that imply the eradication of ticks) with various boundary conditions under the additional
complication of discontinuous diffusion coefficients. These results should be of no surprise to specialists
(they might be even known, but the author did not find a proper reference), because they are related to the
famous Kronig-Penney model (though there the diffusion coefficient is a constant) of quantum mechanics.
These results are tailored in a way to be used for concrete calculations on ticks given below.

The main results are obtained further. They
concern the case of vector-valued diffusions, which allows one to analyse the evolution of the densities of ticks
on various stages of their development. We first suggest some general, though rather rough, sufficient conditions for negativity
by the method of symmetrization. Then we give some more precise conditions, but only in case of two-stage populations
and under rather restrictive assumptions on the control parameters.

In Section \ref{secnumer} we apply our theoretical results to the two examples of ticks populations,
the lone star tick, Amblyoma Americanum, from the south of USA,
and the taiga tick, Ixodes persulcatus, the major type in Russian Siberian forests. There exist quite a lot
of field investigations devoted to these species and thus an abundance of the observations, see e.g. \cite{Mount93}
and \cite{Mount97} and references therein for the American ticks and \cite{Filippova85}, \cite{Korot14} for the taiga tick.
Much less information is available for the European tick, Ixodes ricinus. However, as was many times noted
in the mathematical modelling literature (see e.g. \cite{Gaff07}), there is no consensus on many parameters
required for modelling, and the available parameters
are varying by the types of the territories (upper or bottom land, upland wooded, meadow, etc), humidity and temperature,
seasonal oscillation of hosts, etc. Therefore, numeric estimates of the parameters are obtained
by applying some reasonable averages from the parameters found in the literature.

Section \ref{secproofs} is devoted to the proofs of our results.
In the last section some conclusions are drawn and further perspectives indicated.

\section{Models and theoretical results}
\label{sectheory}

\subsection{Warm-up: one-stage modeling}

In this paper we address the following question: can one place some barriers, organised by some control zones
that are not beneficially for the multiplication of ticks, which would practically stop their propagation
beyond these control zones?

Since we are talking about spatial propagation our model must be necessarily spatially nontrivial. Therefore
we shall apply a diffusion model. As was mentioned above, a diffusion model was developed in  \cite{Caraco02}.
Unlike the model of this paper, our model here will be simpler, in the sense that we will stick to linear
modeling (unlike the reaction-diffusion setting of \cite{Caraco02}), but, on the other hand, it will include
additional complication of being spatially heterogeneous, and precisely the effective control of these heterogeneities
will be the main objective.

\begin{remark} By not including in the model the quadratic (logistic type) terms of the standard nonlinear
model we are not including the competition between species, thus allowing ticks to multiply and propagate under
better conditions than for nonlinear models. Consequently, if we manage to control their propagation in our linear
case (what we are doing here), then it will produce at least the same effect for the corresponding models with
additional nonlinear terms limiting their growth.
\end{remark}

In order to not overcomplicate the story, we will not distinguish infected and noninfected ticks, aiming at their total eradication
(at least locally) and we shall not include explicitly the dynamics of hosts (mice, deers, etc), their overall influence
being reflected just by specifying certain average speed of motion of ticks, which will be effectively modeled as the diffusion
coefficient. Moreover, though the natural state space for the dynamics of ticks is two-dimensional, we shall stick
here to a one-dimensional modeling having in mind the propagation in a certain direction, with the second coordinate
averaged out.

The ticks are known to have several stages of their developments (larvae, nymphs, adult).
To better explain our ideas, we shall first consider the model, where all generations are averaged out, and
next describe the modifications arising when several distinct generations are explicitly included in the model.

Thus we start with the simplest diffusion model for the dynamic of ticks in some area of their habitat
(we denote the derivatives with respect to time and space by a dot and a prime, respectively):
\begin{equation}
\label{eqtickdiff1}
\dot y=ay''+\la y, \quad y(0,x)=y_0(x).
\end{equation}
Here $a>0$ is a diffusion coefficient reflecting the average speed of their random wandering and the coefficient $\la>0$
reflects the speed of their multiplication (usually represented as the difference of average birth rates and average mortality).
The density $y=y(t,x)$ on some interval $x\in [0,R]$ evolves according to \eqref{eqtickdiff1} subject to additional
 boundary conditions, which can be taken as the Dirichlet conditions ($y(t,0)=y(t,R)=0$), or as the Neumann conditions
($y'(t,0)=y'(t,R)=0$) (or more general mixing Robin conditions that we shall not touch here), or rather periodic conditions ($y(t,0)=y(t,R)$, $y'(t,0)=y'(t,R)=0$).

\begin{remark}
\label{remDirper}
 The Dirichlet conditions arise in the situation, when ticks cannot survive at the boundary,
which can be given, for instance, as a water reservoir, or as a temperature barrier.
The Neumann conditions arise in the situation with reflecting barrier,
like a natural or artificial hedge. In other words, this is the situation, when no flux of matter
crosses  the border of the system either from outside or from inside.
Periodic conditions do not seem to have a natural interpretation, but they are often
convenient for the analysis (and therefore are abundantly used in physics), and can be also used as an intermediate step.
In particular, the Dirichlet conditions are clearly the most disastrous for ticks, and therefore, if one can prove the
dying out of ticks for periodic conditions, then the same will hold for the Dirichlet conditions. It is also intuitively
clear (and again widely used in physics) that for large intervals the difference in boundary conditions does not
affect the solutions in any essential way.
\end{remark}

The main point in our modeling is the possibility to attach to the main region $[0,R]$
(the background or the beneficial zone) a control zone of (much smaller length)
$r$, which is not beneficial for ticks, meaning that their reproduction coefficient becomes negative there.
One may also be able to
control the diffusion coefficient shifting it to a different value $b$.
Thus instead of a homogeneous model \eqref{eqtickdiff1} we consider a model with two distinctive regions:
\begin{equation}
\label{eqtickdiff2}
\begin{aligned}
& \dot y_{ben}=ay''_1+\la y_{ben}, \quad x\in (0,R) \\
& \dot y_{nb}=by''_{nb}-\mu y_{nb}, \quad x\in (R,R+r),
\end{aligned}
\end{equation}
where positive $\mu, b, r$ are our control parameters.
To make the problem well posed we should add, as above, some boundary conditions for $x=0$ and $x=R+r$
and moreover, the standard
gluing conditions for diffusions with discontinuous coefficients (see e.g. \cite{SheDec}):
\begin{equation}
\label{eqtickglue}
y_{ben}(R)=y_{nb}(R), \quad ay'_{ben}(R)=by'_{nb}(R).
\end{equation}

Let the operator $L$ be define on the functions of the interval $[0,R+r]$ by the formula
\begin{equation}
\label{eqtickoper1}
Ly(x)=\left\{
\begin{aligned}
& ay''_{ben}(x)+\la y_{ben}(x), \quad x\in (0,R) \\
& by''_{nb}(x)-\mu y_{nb}(x), \quad x\in (R,R+r).
\end{aligned}
\right.
\end{equation}
It is a standard fact from the theory of diffusion operators that,
restricted to the functions satisfying the gluing condition \eqref{eqtickglue} and either Dirichlet,
or Neumann or periodic conditions, the operator $L$ becomes self-adjoint with a discrete
spectrum in the Hilbert space of the square integrable functions $L_2([0,R+r])$ on $[0,R+r]$. Moreover,
this operator is bounded from above, so that there exists a decreasing sequence of eigenvalues
$l_0 \ge l_1 \ge l_2 \ge \cdots$ and the corresponding orthonormal basis $\{\xi_j\}$ in
$L_2([0,R+r])$ (depending on the chosen boundary conditions) such that $L \xi_j=l_j \xi_j$ for all $j$.
This implies that for any initial function $y_0$ from $L_2([0,R+r])$
 the solution to \eqref{eqtickdiff2} with this initial condition
is given by the convergent series
\[
y(x)=\sum_j e^{tl_j}\xi_j(x) \int_0^{R+r} y_0(z) \xi_j(z) dz.
\]
Consequently, if all $\la_j$ are negative, that is the operator $L$ is negative, then this solution
tends to zero for any initial conditions. Thus the condition of eventual eradication of ticks in our model
is the condition of negativity of $L$, or equivalently, the condition of the absence of nonnegative eigenvalues.
In the rare case of vanishing maximal eigenvalue, the solution converges to a finite limit as $t\to \infty$.

\begin{theorem}
\label{th1}
Assume the Dirichlet boundary conditions.
 (i) If
\begin{equation}
 \label{eq1th1}
 \la/a \ge  \pi^2/R^2,
 \end{equation}
$L$ has a positive eigenvalue (independently of $\mu,r$!).

(ii) If
\begin{equation}
 \label{eq2th1}
 \pi^2/(2R)^2 < \la/a <  \pi^2/R^2,
 \end{equation}
the operator $L$ is negative if and only if
 \begin{equation}
 \label{eq3th1}
  -\frac{\tanh (r \sqrt{\mu/b})}{\sqrt{b \mu}}
 > \frac{\tan (R\sqrt{\la/a})}{\sqrt{a\la}}.
\end{equation}

(iii) If
\begin{equation}
 \label{eq4th1}
 \la/a \le  \pi^2/(2R)^2,
 \end{equation}
$L$ has no nonnegative eigenvalues (and is negative in the case of the strict inequality).
\end{theorem}

Statement (i) is the well known (initial) result of the KISS model yielding the critical patch size
\begin{equation}
 \label{criticalpatch}
R_c=\pi\sqrt{a/\la},
\end{equation}
see \cite{KiS53} and \cite{Skellam51}.
Statement (ii) was proved for the limit $r \to \infty$ in \cite{Skellam51}.

Of course for the Neumann and periodic boundary condition the situation is different. There is no critical size,
as the population can survive on any small interval. However, as the following results show, by
introducing appropriate control zones, one can eradicate the population in a similar way to the above.

\begin{theorem}
\label{th1a}
Assume the Neumann boundary conditions.
 (i) If
\begin{equation}
 \label{eqNeu1}
 \la/a \ge  \pi^2/(2R)^2,
 \end{equation}
$L$ has a positive eigenvalue (independently of $\mu,r$!).

(ii)  Otherwise, $L$ is negative if and only if
 \begin{equation}
 \label{eqNeu2}
 \sqrt{\mu b} \tanh (r \sqrt{\mu/b}) > \sqrt{\la a}  \tan (R\sqrt{\la/a}),
\end{equation}
so that one can achieve negativity by choosing appropriate $\mu$ and $r$.
\end{theorem}

\begin{theorem}
\label{th2}
Assume the periodic boundary conditions.
 (i) If \eqref{eq1th1} holds, then
$L$ has a positive eigenvalue.
(ii) Otherwise, $L$ is negative if and only if
 \begin{equation}
 \label{eq1th2}
 \sqrt{\mu b} \tanh (\frac{r}{2} \sqrt{\mu/b}) > \sqrt{\la a}  \tan (\frac{R}{2}\sqrt{\la/a}).
\end{equation}
\end{theorem}

One sees that  Theorem \ref{th2} provides the same condition as Theorem \ref{th1a} though with $R/2$ and $r/2$
instead of $R$ and $r$. On the other hand, the boundary KISS value of \eqref{eq1th1} is the same for the Dirichlet
and the periodic conditions. Moreover, one can check explicitly from conditions \eqref{eq1th2} and \eqref{eq3th1}
that, as expected (see Remark \ref{remDirper}), if the operator $L$ is
negative under the periodic condition, then it is also negative  under the Dirichlet conditions.

 Theorems \ref{th1} and \ref{th2} provide exact information on how $\mu,b,r$ can be tuned in order
 to achieve negativity of $L$ and hence the eventual eradication of ticks. We see, in particular, that if
$R$ does not exceed certain critical level, one can eradicate the ticks by introducing a control zone of
arbitrary small length $r$,
if a sufficiently high level of mortality $\mu$ can be imposed on this control zone.

On the other hand, the theorems show that if $\la R^2/a$ is large enough, then no control zone can
 efficiently influence the global growth of ticks population. This observation leads to a natural idea that,
 in case of a large territory, one can fight with the growth of ticks by arranging several small
 control zones placed in a periodic (patched) fashion.
 As we are going to show, if this is organised in a way that each pair of adjacent zones satisfies the
 conditions of the previous theorem, the ticks will be eradicated in the whole patched territory.

Thus assume that our territory of habitat is represented as the interval $[0,K(R+r)]$,
with some natural $K$, and that it is decomposed into $2K$ subintervals, odd
$I_1=[0,R]$, $I_3=[R+r, 2R+r]$, $\cdots$, $I_{2K-1}=[(K-1)(R+r), KR+(K-1)r]$, so that
$I_{2k-1}=[(k-1)(R+r), kR+(k-1)r]$ for $k=1, \cdots, K$,
 and even $I_2=[R, R+r], \cdots , I_{2K}=[KR+(K-1)r, K(R+r)]$, so that
$I_{2k}=[kR+(k-1)r, k(R+r)]$ for $k=1, \cdots, K$. Assume that on the long odd intervals
we have some background parameters  $\la>0, a>0$, and on the short even intervals we have some
(controllable) parameters $b>0,\mu>0$, so that the diffusion is given by the first and second
equations of \eqref{eqtickoper1} on even  and odd intervals, respectively. Of course the gluing
conditions \eqref{eqtickglue} are supposed to hold on the interface of all intervals.
Theorem \ref{th2} (valid for the case $K=1$) can now be extended to the case of arbitrary $K$
as follows.

\begin{theorem}
\label{th3}
Assume the periodic boundary conditions for such diffusion on $[0,K(R+r)]$. Then for any $K$ the
conditions of Theorem \ref{th2} provide also the conditions for the corresponding diffusion operator
$L$ on $[0,K(R+r)]$ to be negative or not.
\end{theorem}

\begin{remark}
As we already mentioned, if all positive solutions for the periodic conditions tend to zero, as $t\to \infty$,
then the same holds for the Dirichlet conditions. It is intuitively clear, as the Dirichlet conditions are less
beneficial for survival. Formally it follows from the representations of solutions in terms of the Feynmann-Kac
formula.
\end{remark}

\subsection{Several-stage modeling: trivial case}

So far we have considered the situation with all stages of ticks' lives averaged out.
In more precise modeling one has to take into account the presence of several stages.
For ticks these are eggs, larvae, nymphs and adults.
With some reasonable averaging one can reduce
the consideration of the lifespan of ticks to the two basic periods, from eggs to nymphs,
and from nymphs to hatching female adults. On the other hand, a more detail analysis, can include
not only the stages, but their time developments. Namely, say, nymphs can develop to adults in the
same season as their own molting takes place, or after a diapause (wintering), so one can distinguish
not only the stage, but also whether it develops in one year or two years.

Let us consider the general case of $n$
stages (could be also generations for other species). The basic equation \eqref{eqtickdiff1}
is generalised to the vector-valued equation
\begin{equation}
\label{eqtickdiffvec1}
\dot y=ay''+My, \quad y(0,x)=y_0(x),
\end{equation}
where $y\in \R^n$ and $M$ is the birth-and-death matrix showing the progression of
ticks from their birth through
their various stages of development. This matrix $M=M_n$ has the standard form
(used by many authors, see e.g.  \cite{Mwa00} or \cite{Zhao12}) reflecting the sequential propagation
through various stages. They have elements
\begin{equation}
\label{eqmatrixstages}
M_{jj}=-m_j, \, j=1, \cdots, n, \quad M_{j+1,j}=b_j, \, j=1, \cdots , n-1, \quad M_{1n}=b_n,
\end{equation}
with all other elements vanishing, where all parameters $\mu_j$ and $b_j$ are positive.

For instance, for dimensions $2$ and $3$, these matrices have the form
\begin{equation}
\label{eqmatbirthdeath}
M_2=
\left(
\begin{aligned}
& -m_1 \quad  b_2 \\
& \,\,\, b_1 \quad -m_2
\end{aligned}
\right),
\quad
M_3=
\left(
\begin{aligned}
& -m_1 \quad \,\,\, 0 \quad \quad b_3 \\
& \quad b_1 \quad -m_2 \quad \,\, 0 \\
& \quad 0 \quad \quad \,\, b_2 \quad -m_3
\end{aligned}
\right).
\end{equation}

It is easy to see that a matrix $M$ of type \eqref{eqmatrixstages} always has a real eigenvalue.

Let us assume now that we can set a control zone with the increased death rates for ticks. Namely,
let us consider the extension of equation \eqref{eqtickdiff2} of the form

\begin{equation}
\label{eqtickdiffvec2}
\begin{aligned}
& \dot y_{ben}=ay''_{ben}+My_{ben}, \quad x\in (0,R) \\
& \dot y_{nb}=by''_{nb}+(M-\mu \1)y_{nb}, \quad x\in (R,R+r),
\end{aligned}
\end{equation}
where $\1$ is the unit matrix. Notice that $\mu $ here has a slightly different meaning as in \eqref{eqtickdiff2},
referring to the relative decrease in the death rates.
The corresponding diffusion operator becomes matrix-valued:
\[
Ly(x)=[ay''(x)+My(x)]\1_{[0,R]}(x)+[by''(x)+(M-\mu \1)]\1_{[R,R+r]}(x).
\]

It turns out that the results of Theorems \ref{th1}-\ref{th3} have a straightforward extension
 to this case with the role of $\la$ played by the largest eigenvalue of $M$.

 \begin{theorem}
\label{thvec1}
Assume that the maximal real eigenvalue $\La_1$ of the matrix $M$ in \eqref{eqtickdiffvec2}
is positive (if all eigenvalues of $M$ have negative real parts,
then the population would die out even in the background territory,
the case of no interest to us) and all other eigenvalues $\La_j$ have negative
real part and are different (the latter conditions are technical simplifications
that are not essential). Assume also that $\mu>\La_1$ (otherwise the ticks
could survive even in the control zone alone).

Assume the Dirichlet boundary conditions for $L$.
 (i) If
\begin{equation}
 \label{eq1thvec1}
 \La_1/a \ge  \pi^2/R^2,
 \end{equation}
$L$ has a positive eigenvalue.
(ii) If
\begin{equation}
 \label{eq2thvec1}
\pi^2/(2R)^2 < \La_1/a < \pi^2/R^2,
 \end{equation}
the operator $L$ is negative if and only if
 \begin{equation}
 \label{eq3thvec1}
  -\frac{\tanh (r \sqrt{(\mu-\La_1)/b})}{\sqrt{b (\mu-\La_1)}}
 > \frac{\tan (R\sqrt{\La_1/a})}{\sqrt{a\La_1}}.
\end{equation}

(iii) If
\begin{equation}
 \label{eq4thvec1}
 \La_1/a \le  \pi^2/(2R)^2,
 \end{equation}
$L$ has no nonnegative eigenvalues (and is negative in the case of the strict inequality).
\end{theorem}

Assume now that the territory of habitat is represented by the interval $[0,K(R+r)]$
(with some natural $K$), decomposed into $2K$ subintervals in the same way as formulated
before Theorem \ref{th3}, and that on the odd and even intervals our diffusion follows the first
and the second equation of \eqref{eqtickdiffvec2}, respectively, with the usual gluing condition
on the interfaces.

\begin{theorem}
\label{thvec2}
Let the conditions of Theorem \ref{thvec1} for $M$ hold and $K$ be arbitrary.
Assume that the periodic boundary conditions are chosen.
 (i) If \eqref{eq1thvec1} holds, then
$L$ has a positive eigenvalue.
(ii) Otherwise the operator $L$ is negative if and only if
 \begin{equation}
 \label{eq1thvec2}
 \sqrt{(\mu-\La_1) b} \tanh (\frac{r}{2} \sqrt{(\mu-\La_1)/b}) > \sqrt{\La_1 a}  \tan (\frac{R}{2}\sqrt{\La_1/a}).
\end{equation}
\end{theorem}

\subsection{Several-stage modeling: advanced case}

This section contains our main theoretical results.

So far we have looked at the case when diffusion coefficients differ in beneficial and control zones,
but are independent of the stage. Of course, usual averaging allows one to apply this model for ticks.
However, for many types of ticks, larvae and nymphs use small rodents as hosts, while adult ticks use
large mammals, like deers, or birds, so that the displacement, and hence the diffusion coefficient differ drastically
for the adults and the earlier stages of ticks. Hence it is more natural to use the model with different
diffusions on different stages. Moreover, the matrices $M$ specifying the birth and death rates,
can be of course quite different for the background and control zones, thus differing not only by a multiple
of the unit matrix, as in \eqref{eqtickdiffvec2}.

Let us start with a simple extension of the KISS model. Namely,
let us consider the extension of equation \eqref{eqtickdiffvec1} with variable diffusion,
that is, the equation
\begin{equation}
\label{eqtickdiffvec1n}
\dot y=Ay''+My, \quad y(0,x)=y_0(x),
\end{equation}
where $A$ is a diagonal matrix with diagonal elements $a_j>0$, $j=1, \cdots, n$.

\begin{theorem}
\label{thvec3}
Let $M$ be a matrix of type \eqref{eqmatrixstages}. Suppose that the maximal eigenvalue
$\La_1=\La(M,A)$ of the matrix $A^{-1}M$ is strictly positive.
Then the critical patch size equals
$R_c=\pi / \sqrt{\La_1}$. That is, there exists a nontrivial solution to the equation $Ay''+My=Ey$
with some positive $E$ if and only if $R>R_c$.
\end{theorem}

The story becomes more complicated when we put together the original and a control zones. It seems difficult
to expect here explicit necessary and sufficient conditions for the negativity of the spectrum, like in the cases,
analyzed above. However, reasonable sufficient conditions can be obtained, as we are going to show now.

Let us consider the system
\begin{equation}
\label{eqtwopart}
 \begin{aligned}
 & \dot y_{ben}=A_{ben}y''_{ben} +M_{ben}y_{ben}, \quad y\in (-R/2,R/2), \\
 & \dot y_{nb}=A_{nb}y''_{nb} +M_{nb}y_{nb}, \quad y\in (R/2, r+R/2),
\end{aligned}
\end{equation}
with symmetric positive matrices $A_{ben}$, $A_{nb}$ and arbitrary matrices $M_{ben}$, $M_{nb}$,
where matrices $M_{nb}$, $A_{nb}$ are supposed to depend on some control parameters.
We assume the periodic boundary conditions and the usual gluing
conditions ($y_{nb}=y_{ben}$, $A_{nb}y'_{nb}=A_{ben}y'_{ben}$) on the interface.

Equations for the eigenvalues are
 \begin{equation}
\label{eqtwopart1}
 \begin{aligned}
 & A_{ben}y''_{ben} +M_{ben}y_{ben}=Ey_{ben} \\
 & A_{nb}y''_{nb} +M_{nb}y_{nb}=Ey_{nb}.
\end{aligned}
\end{equation}

\begin{theorem}
\label{thvec4}
Let exactly $k$ out of $n$ eigenvalues $\la_1 \ge \cdots \ge \la_n$ of the symmetric matrix
$N_{ben}=(M_{ben}A_{ben}^{-1}+A_{ben}^{-1}M_{ben}^T)/2$ are positive, and all
eigenvalues $\mu_1 \ge \cdots \ge \mu_n$ of the symmetric matrix
$N_{nb}=(M_{nb}A_{nb}^{-1}+A_{nb}^{-1}M_{nb}^T)/2$ are negative (by $T$ we denote the transposition). If
\begin{equation}
\label{eq1thvec4a}
R\le R_c^{sym}=\pi/\sqrt{\la_1},
\end{equation}
and
\begin{equation}
\label{eq1thvec4}
 \frac{\sqrt{|\mu_1|} \sinh (r\sqrt{|\mu_1|})}
{1+\cosh(r\sqrt{|\mu_1|})}
>\frac{2Rnk \la_1}{1+\cos (R\sqrt{\la_1}))},
\end{equation}
then system \eqref{eqtwopart1} has no solutions with non-negative $E$.
\end{theorem}

This result extends automatically to the case of the territory represented by the interval $[0,K(R+r)]$,
with some natural $K$, which is decomposed into $2K$ subintervals in the same way as formulated before
Theorem \ref{th3} (see also Theorem \ref{thvec2}),
so that on the odd and even intervals the diffusion follows the first
and the second equation of \eqref{eqtwopart}.

Notice that the r.h.s. of \eqref{eq1thvec4} increases as $\sqrt{|\mu_1|}$ like in our previous results, but
the l.h.s. is of order $\la_1$ (for $R$ far away from $R_c^{sym}$), unlike $\sqrt{\la_1}$ in Theorem \ref{th2}.

In Theorem \ref{thvec3} the condition of the existence of a positive eigenvalue is given in terms
of the maximal eigenvalue
$\La_1$ of the matrix $A^{-1}M$ and in Theorem \ref{thvec4} a sufficient condition for negativity is linked with
the maximal eigenvalue of the corresponding symmetrised matrix. Hence from the combination of these two theorems
nothing can be said for sizes $R$ such that
\begin{equation}
\label{eq1thvec4b}
\pi/\sqrt{\la_1}=R_c^{sym} < R < R_c=\pi/\sqrt{\La_1}.
\end{equation}

In order to work with these sizes the effective method of symmetrization used in  Theorem \ref{thvec4} cannot be applied,
and the analysis in arbitrary dimension seems to be quite involved. Working in dimension $n=2$, where explicit formulas for
 eigenvectors allow for a rather detailed analysis, possible results seem to depend strongly on the structure of diffusion
 coefficients and the matrix $M$. We present one such result valid for a range of coefficients that can be applied
 to tick populations.

Let us look at the equations
\begin{equation}
\label{eq2dimvar}
 \begin{aligned}
 & \dot y_{ben}=Ay''_{ben} +M_{ben}y_{ben} \\
 & \dot y_{nb}=aAy''_{nb} +M_{nb}y_{nb}
\end{aligned}
\end{equation}
and the corresponding eigenvalue problem
\begin{equation}
\label{eq2dimvarst}
 \begin{aligned}
 & Ay''_{ben} +M_{ben}y_{ben}=Ey_{ben} \\
 & aAy''_{nb} +M_{nb}y_{nb}=Ey_{nb},
\end{aligned}
\end{equation}
with the usual gluing condition: $y_{ben}$ coincides with $y_{nb}$ and $y'_{ben}$ coincides with $ay'_{nb}$
on the interface of two regions.

Let $n=2$.
For a $E\ge 0$, denote by $\La_1(E)\ge\La_2(E)$  the eigenvalues of the matrix $N_{ben}(E)=A^{-1}M_{ben}-A^{-1}E$,
and by $\mu_1(E)\ge \mu_2(E)$ the eigenvalues of the matrix
$N_{nb}(E)=(A^{-1}M_{nb}-A^{-1}E)/a$. Let $\La_1=\La_1(0)>0>\La_2(0)$ and
 $E_0>0$ be such that $\La_1(E_0)=0$. Assume that
$\La_1(E)\ge 0>\La_2(E)$ and $0>\mu_1(E)>\mu_2(E)$ for all $E\in [0,E_0]$.
Let $\{v_j(E)\}$ and $\{w_j(E)\}$ be the corresponding bases of eigenvectors
of $N_{ben}(E)$ and $N_{nb}(E)$, and
$c(E)=(c_{ij}(E))$ be the matrix that takes the basis $\{v_j(E)\}$ to the basis $\{w_j(E)\}$.

\begin{theorem}
\label{thvec5}
Assume the bases of eigenvectors can be chosen in such a way that
 \begin{equation}
\label{eq2thvec5}
c_{12}(E)c_{21}(E)\le 0, \quad c_{11}(E)c_{22}(E)\ge 0
\end{equation}
for $E\in [0,E_0]$. Then, if
 \begin{equation}
\label{eq1thvec5}
a\sqrt {\mu_1(0)} \tanh (\sqrt{\mu_1(0)}\frac{r}{2})> \sqrt \La_1 \tan (\sqrt{\La_1}\frac{R}{2}),
\end{equation}
there are no positive solutions to equation \eqref{eq2dimvarst}.
\end{theorem}

 Comparing \eqref{eq1thvec5} with \eqref{eq1th2} we see that under \eqref{eq2thvec5}
our two-dimensional condition of negativity is the exact extension of the one-dimensional case.

\begin{remark}
The technically convenient assumptions \eqref{eq2thvec5}
(and especially the first of these two) are not very natural
and should not hold in a general situation. However, as we show below, they can be achieved
by an appropriate proportional change of the elements of the birth-and-death matrix.
\end{remark}

 Let us now look at a concrete situation when  \eqref{eq2thvec5} holds. Namely, consider the system
\eqref{eq2dimvar} where $a=1$, the matrix $A$ is diagonal with the diagonal elements $a_1,a_2$
$M_{ben}$ has the form $M_2$ from \eqref{eqmatbirthdeath} and
$M_{nb}$ has the same form with $\tilde m_j$ and $\tilde b_j$ instead of $m_j$ and $b_j$.
Thus
\[
N_{ben}(0)= A^{-1}M_{ben}=
\left(
\begin{aligned}
& -a_1^{-1}m_1 \quad  a_1^{-1} b_2 \\
& \,\,\, a_2^{-1}b_1 \quad -a_2^{-1}m_2
\end{aligned}
\right), \quad
N_{nb}(0)= A^{-1}M_{nb}=
\left(
\begin{aligned}
& -a_1^{-1} \tilde m_1 \quad  a_1^{-1} \tilde b_2 \\
& \,\,\, a_2^{-1}\tilde b_1 \quad -a_2^{-1} \tilde m_2
\end{aligned}
\right).
\]

Assume that we can control the non-beneficial zone by a proportional decrease of
the reproduction coefficients, that is calibrating
 \begin{equation}
\label{eq0thvec6}
\tilde b_j=\om b_j, \quad j=1,2,
\end{equation}
by choosing an appropriate parameter $\om\in(0,1)$ and choosing $\tilde m_j \ge m_j$ in such a way
that
 \begin{equation}
\label{eq00thvec6}
\tilde m_1-\tilde m_2\ge m_1-m_2.
\end{equation}
For real ticks, $a_1$ is usually  much less than $a_2$ (see numeric examples below).
Hence the assumptions of the next theorem are fully relevant.

\begin{theorem}
\label{thvec6}
Assume that $\det N_{ben}(0)=m_1m_2-b_1b_2<0$ (this enures that $\La_1(0)>0>\La_2(0)$) and
 \begin{equation}
\label{eq1thvec6}
a_1^{-1} \ge a_2^{-1}, \quad   a_1^{-1} m_1 \ge  a_2^{-1} m_2.
\end{equation}
Assume that $\tilde b_j$ and $\tilde m_j$ satisfy \eqref{eq0thvec6} and \eqref{eq00thvec6}.
Then the conditions of Theorem \ref{thvec5} (namely, equation \eqref{eq2thvec5}) hold, so that
\eqref{eq1thvec5} is sufficient for the negativity
 of all eigenvalues (and hence for the eradication of ticks).
\end{theorem}

\section{Numerical results with real life data}
\label{secnumer}

\subsection{One-stage modeling with the North American ticks}

In US, the predominant types of ticks are the lone star tick (Amblyoma Americanum) in the south,
the blacklegged tick (Ixodes scapularis, formerly called the deer tick) and the Americal dog tick
 (Dermacenter Variabilis) in the north. One can find lots of experimental
 research on the various parameters needed for modeling, which depend on many factors. Here, to estimate
 the birth and death rates, we will employ the averaged parameters for the lone star tick from
 \cite{Gaff07}. Namely, the time unit is a month. Taking 2000 as the average number of eggs
 (occurring once in two year), $70\%$ survival rate and sex ratio $1:1$, we get the birth rate
of female larvae per month to be $0.5\times 0.7 \times 2000/24$. With the out of host survival
average 0.85 and the probability to find the host 0.03 (recall that without finding a host for
a full meal blood no further developments of ticks is possible)
the total survival rate (of females) per month is calculated in \cite{Gaff07} as
\[
0.5\times 0.7 \times 2000/24 \times 0.03 \times 0.85=0.75.
\]
Death rate is estimated as $0.01$ in woods and as $0.1$ in grass. Sticking to the case of grass, we get the parameter $\la$ in
equation \eqref{eqtickdiff1} to be $\la=0.75-0.1=0.65$ (measured in month${}^{-1}$).

To estimate the diffusion coefficient $a$ is a more difficult task.

We shall use the standard method, used both in physics and ecology (see e.g. \cite{Holmes93}),
where the mean squared displacement $M^2(t)$ during a time $t$ is estimated as
\begin{equation}
\label{diffcoeff}
M^2(t)=2a t.
\end{equation}
The life time of the majority of
lone star ticks is known to be 2 years, but sometimes it is completed in one year.
Thus we can assume approximately that it has two rides per year.
 The main hosts for the lone star ticks (and many other American Ticks) are the white tailed deers
 that could travel for many kilometers during 3 - 4 days needed for a tick to get its blood meal.
 One can choose about 10 km as the reasonable estimate for an average distance per
 a ride (see \cite{Ostfeld}). Note that the displacement due to ticks' own movement is negligibly small
 compared to their displacement by the hosts.

 Then the total squared displacement is $M^2(12)=(2\times 10)^2=400$, and
 therefore $a=400:24\approx 16.67$, measured in km$^2$ per month.  Thus for the critical patch size we get
 \begin{equation}
\label{critpatchamer}
 R_c=\pi \sqrt{16.67/0.65} \approx 15.9 \, \text{(km)}.
 \end{equation}
 Theorems \ref{th1} - \ref{th3} can be used to define the exact relation between the length
 of a control zones and the death rate  in it, which is needed for the eradication of ticks.

 For instance, let us apply Theorem \ref{th3}.

 Let us choose $R=14$ (km) (which is reasonably close to the critical size $R_c$)
 and the length of the control zone $r=1$ (km), and let $b=a$ for simplicity.
 Then condition \eqref{eq1th2} for the eradication of ticks
 gets the following numeric form
  \begin{equation}
 \label{critpatchamer1}
 \sqrt{16.67 \mu } \tanh (0.5 \sqrt{\mu/16.67}) > 17.03,
\end{equation}
so that the minimal required death rate $\mu$ is rather high, about 1958.

\subsection{One-stage modeling with the taiga ticks}

Let us exploit here another time unit, choosing it to be one year.
The life cycle of the taiga tick varies from 3 to 6, so that 4 years can be taken as
an approximate average. Recall that in order to have a molting and to turn from one
stage to another (larvae to a nymph, nymph to an adult, and finally to hatch)
a tick must have a ride on a host with a full meal. Hence with 3- 4 year cycle
 a tick makes on average a single ride per year.

Taking, as above, 10 km as an average distance per a meal (and thus per year) we get
 from \eqref{diffcoeff} that $a=100/2=50$, measured in km$^2$ per year.

As the mean survival rate we choose $\la=2$, measured in year$^{-1}$, which is the
approximate value used in \cite{Vshiv13}, where it was shown to produce a reasonable fit
to the experimental data available.  Thus for the  critical patch size we get
 \begin{equation}
\label{critpatchtaiga}
 R_c=\pi \sqrt{50/2} \approx 15.7 \, \text{(km)},
 \end{equation}
 giving approximately the same result (!), as in the calculations above, based on the American experimental data.


\subsection{Two-stage modeling with the taiga ticks}

One of the ways to estimate the average travel distance of ticks on hosts can be obtained
by using chemical treatments of a controlled territory and looking for how far ticks can
penetrate into the treated territory from the uncontrolled zone. Theses studies indicate
(see Section VIII.3 of \cite{Filippova85}) that, for a taiga tick, one can choose about
10 km as the average displacement for adult ticks during their meals (supporting the number
used above and taken from the literature on American ticks) and about 1.5 km for nymphs and larvae.
Therefore the diffusion coefficient $a=100/2=50$, used above for the whole population of ticks, in a more detailed analysis
becomes the diffusion coefficient $a_2=50$ for the adult ticks only, while for all lower stages it can be estimated as only
$a_1=1.5^2/2\approx 1.1$.

Let us work with the two-stage model, where the first stage represent the species grown up to the well-fed nymph,
and the second stage represent adults. We shall use here the data for the taiga tick from the Sayans mountains of Siberia,
Krasnoyarsk region,
as presented in  the very detailed observations of \cite{Korot14}. These observations show that only about $2.4\%$ of
female tick get their meal and that about $1/2$ of well-fed females produce eggs in the amount of about 5 thousand each.
Taking also the standard sex ratio as $1:1$, it follows that the potential number for the next generation is
\[
p=5000 \times 0.024 \times 0.5^2=30
\]
 per an adult tick. From this amount about $43\%$ survive producing hungry nymphs in the amount
 of $0.43p=12.9$. Only $8.2\%$ of the potential $p$ survive to the stage of an adult well-fed nymphs
 (because of the high death rate on this stage), that is in total remain $0.082p=2.46$ per an original adult tick.

 Next, the death rate from a well-fed nymph to an adult tick is about $14.6\%$, and the
 deathrate during a winter is about $39\%$. Thus from a well-fed nymph one can expect about
 \[
 0.854 \times 0.61\approx 0.52
 \]
 hungry adults to appear next spring. Thus the birth-and-death matrix $M$ and
 the diffusion matrix $A$ in \eqref{eqtickdiffvec1n} take the concrete form
 \begin{equation}
\label{eqmatrkorot}
M= \left(
\begin{aligned}
 & -1 \quad 2.46 \\
 & 0.52 \quad -1
\end{aligned}
\right), \quad
A= \left(
\begin{aligned}
 & 1.1 \quad 0 \\
 & \,\, 0 \quad 50
\end{aligned}
\right).
\end{equation}

Therefore
\[
A^{-1}M=\left(
\begin{aligned}
 & -0.91 \quad 2.46\times 0.91 \\
 & 0.52\times 0.02 \quad -0.02
\end{aligned}
\right)
=\left(
\begin{aligned}
 & -0.91 \quad 2.24 \\
 & \,\,\,  0.01 \quad -0.02
\end{aligned}
\right)
\]
For the positive eigenvalue $\La_1$ of this matrix we get approximately that $\sqrt{\La_1}=0.067$
and thus the critical patch size of Theorem \ref{thvec3} equals
 \begin{equation}
\label{critpatchtaiga2}
R_c=\pi/0.067 \approx 46.9.
\end{equation}
This value is essentially larger than the values \eqref{critpatchamer} and  \eqref{critpatchtaiga}, calculated for
one-stage models and data (and taken from different sources). This can be expected. In fact, more detailed
observations and calculations based on the regions of approximate  equilibrium should reflect only slight possible average
growth rates and hence weaker requirements for control zones. Theorem \ref{thvec6} can be used to assess the required decrease
in the reproduction rates of control zones that would ensure the eradication.

Choosing, say, $R=40,r=1$,
the condition of negativity \eqref{eq1thvec5} gets the following numeric form:

 \begin{equation}
\label{critpatchtaiga2a}
\sqrt {\mu_1(0)} \tanh (\sqrt{\mu_1(0)}/2)> 0.067 \tan (1.34)=0.29.
\end{equation}
Therefore $|\mu_1|$ must be larger than about $0.64$. This condition yields the
corresponding estimates for the required $\tilde m_j, \tilde b_j$.

On the other hand,
\[
\frac12 (A^{-1}M+M^T A^{-1})
=\left(
\begin{aligned}
 & -0.91 \quad 1.125 \\
 & \,\,\,  1.125 \quad -0.02
\end{aligned}
\right).
\]
Consequently, one gets for the highest eigenvalue $\la_1$ of this symmetrised matrix
that $\sqrt \la_1 \approx 0.86$. Therefore
 \begin{equation}
\label{critpatchtaiga3}
R_c^{sym}=\pi/0.863 \approx 3.64.
\end{equation}

This value is much less than values \eqref{critpatchamer} and  \eqref{critpatchtaiga}. This corroborates the idea
that the results of Theorem \ref{thvec4}, which are very natural theoretically, must be used cautiously in
situations with highly nonsymmetric birth-and-death matrices,
as the symmetrization leads to strong distortions for such data.

\begin{remark}
\label{remimprovesym}
One of the ways to enhance the practical application of Theorem \ref{thvec4}
might be the introduction of additional life
stages (increasing dimension) that would make the birth-and-death matrix more symmetric.
\end{remark}

\section{Proofs}
\label{secproofs}

{\it Proof of Theorem \ref{th1}}.

 Negativity of the operator $L$ means that $L$ has no non-negative eigenvalues.
If such an eigenvalue $E$ exists, then it solves the stationary problem $Ly=Ey$.
Since the operator $y\to y''$ with
Dirichlet boundary conditions is negative, it follows that $E\in[0,\la)$. Then the
general solutions in the two domains are
\begin{equation}
\label{eqgensolbas}
y_{ben}=C\sin \left(\sqrt{\frac{\la-E}{a}}x\right)+D \cos \left(\sqrt{\frac{\la-E}{a}}x\right), \quad x\in (0,R),
\end{equation}
\begin{equation}
\label{eqgensolcon}
y_{nb}=A \exp\{ \sqrt{\frac{\mu+E}{b}}x\}+B \exp\{ - \sqrt{\frac{\mu+E}{b}}x\}, \quad x\in (R,R+r).
\end{equation}

The Dirichlet boundary condition $y_{ben}(0)=0$ implies that $D=0$. Since eigenfunctions are defined
up to a multiplicative constant, we can set $C=1$, so that
\[
y_{ben}=\sin \left(\sqrt{\frac{\la-E}{a}}x\right).
\]
The Dirichlet boundary condition $y_{nb}(R+r)=0$ yields
\[
B=-A \exp\{ 2\sqrt{\frac{\mu+E}{b}}(R+r)\}.
\]

By the gluing conditions $y_{ben}(R)=y_{nb}(R)$ and $ ay'_{ben}(R)=by'_{nb}(R)$ on the interface,
\[
\sin \left(\sqrt{\frac{\la-E}{a}}R\right)=A \exp\{ \sqrt{\frac{\mu+E}{b}}R\}+B \exp\{ - \sqrt{\frac{\mu+E}{b}}R\},
\]
\[
\sqrt{a(\la-E)} \cos \left(\sqrt{\frac{\la-E}{a}}x\right)
=A \sqrt{b(\mu+E)} \exp\{ \sqrt{\mu+E}R\}-B \sqrt{b(\mu+E)}\exp\{ -\sqrt{\mu+E}R\}.
\]

Thus
\[
\sin \left(\sqrt{\frac{\la-E}{a}}R\right)
=A \left[\exp\{ \sqrt{\frac{\mu+E}{b}}R\}- \exp\{ \sqrt{\frac{\mu+E}{b}}(R+2r)\}\right],
\]
\[
\sqrt{a(\la-E)} \cos \left(\sqrt{\frac{\la-E}{a}}x\right)
\]
\[
=A \sqrt{b(\mu+E)}
\left[\exp\{ \sqrt{\frac{\mu+E}{b}}R\}+ \exp\{ \sqrt{\frac{\mu+E}{b}}(R+2r)\}\right],
\]
or
\[
-\frac{\tanh \left(r \sqrt{\frac{\mu+E}{b}}\right)}{\sqrt{b(\mu+E)}}
= \frac{\tan \left(\sqrt{\frac{\la-E}{a}}R\right)}{\sqrt{a(\la-E)} }.
\]
 Changing to $x=(\la-E)/a\in (0, \la/a]$ yields
  \begin{equation}
 \label{eq1di}
 -\frac{\tanh (r \sqrt{(\la+\mu-ax)/b})}{\sqrt{b(\la+\mu-ax)}} = \frac{\tan (R\sqrt{x})}{a\sqrt{x} }.
\end{equation}

The l.h.s. is negative decreasing on $[0, \la/a]$. Thus, if
 $\la/a \le \pi^2/(4R^2)$,
 there is no solutions $x\in (0, \la/a]$, because the r.h.s. is positive on this interval.

Moreover, the l.h.s. is decreasing on $[0, \la/a]$ from
\[
 -\frac{\tanh (r \sqrt{(\la+\mu)/b})}{\sqrt{b(\la+\mu)}} \,\, \text{to} \,\,
 -\frac{\tanh (r \sqrt{\mu/b})}{\sqrt{b \mu}}.
 \]
 Thus, if $\la/a \ge \pi^2/R^2$,
 there is a solution $x\in (0, \la/a]$, because the l.h.s. necessarily intersects with
 the part of the r.h.s. increasing from $-\infty$ to $0$. Finally, if
 \[
 \frac{\pi^2}{4R^2} < \la/a <  \frac{\pi^2}{R^2},
 \]
 a solution exists if and only if the value of the l.h.s of \eqref{eq1di} at
 $x=\la/a$ lies below the graph of the r.h.s. yielding condition \eqref{eq3th1}.

\vspace{3mm}
{\it Proof of Theorem \ref{th1a}}. Similar calculations to Theorem \ref{th1} show that
the condition for $E\in [0,\la)$ to be an eigenvalue writes down as the equation

\begin{equation}
 \label{eqNeu0pr}
\sqrt{a(\la-E)} \tan  \left(\sqrt{\frac{\la-E}{a}} R\right)
=\sqrt{b(\mu+E)} \tanh \left(\sqrt{\frac{\mu+E}{b}} r\right),
\end{equation}
or, in terms of $x=(\la-E)/a$, as
 \begin{equation}
 \label{eqNeu1pr}
 \sqrt{(\mu+\la-ax)b} \tanh \left(r \sqrt{\frac{\mu+\la-ax}{b}}\right) = a\sqrt{x}  \tan (R\sqrt{x}).
\end{equation}

The r.h.s. is positive increasing from $0$ to $\infty$ on $[0, \pi^2/4R^2)$.
The l.h.s. is positive decreasing on $[0,\la]$ between two positive values (or between a positive value and $0$
if $\mu=0$). Hence if
\[
\la/a \ge \pi^2/(2R)^2,
\]
there is a positive eigenvalue for any $\mu, r$. Otherwise, a positive eigenvalue does not exist if and only
\eqref{eqNeu2} holds.

\vspace{3mm}
{\it Proof of Theorem \ref{th2}}.
We again look at the system of equations
\[
\dot y_{ben}=ay''_{ben}+\la y_{ben}, \quad \dot y_{nb}=b y''_{nb}-\mu y_{nb}.
\]
It is convenient to choose coordinates so that the first equation holds for $x\in (-R/2,R/2)$ and the second for
$x\in (R/2,R/2+r)$. Here the key point for the analysis is the symmetry. Namely, our system writes down as the
equation $\dot y=Ly$ on the torus $y\in [-(R+r)/2, (R+r)/2]$ with the identified right and left points. Alternatively
one can think of $y$ as periodic functions (with period $R+r$) on the whole $\R$.
 Here $L$ is the linear operator such that $L\RC=\RC L$, where $\RC$ is the reflection
operator: $\RC f(x)=f(-x)$. Hence if $y$ is an eigenfunction of $L$, then so is also the function $\RC y$.
  But $L$ is an operator of the Schr\"odinger type, and it is well known from quantum mechanics
(and actually very easy to prove, see e.g. \cite{LaLiQM}) that the eigenvalues of one-dimensional Scr\"odinger operators are
non-degenerate, that is, there may exist only one eigenfunction to each eigenvalue, up to a multiplier. Hence
$\RC y=\pm y$. But the minus sign would contradict the
continuity and thus $\RC y=y$. Thus an eigenfunction of $L$, which is the solution to the equation $LY=Ey$
must be an even function. Similarly it has to be symmetric under the reflection with respect to the point $(R+r)/2$.

The equation $Ly=Ey$ writes down more explicitly as the system
\[
ay''_{ben}+\la y_{ben}=Ey_{ben}, \quad b y''_{nb}-\mu y_{nb}=Ey_{nb},
\]
or equivalently
\[
y''_{ben}+\frac{\la-E}{a} y_{ben}=0, \quad y''_{nb}-\frac{\mu +E}{b}y_{nb}=0.
\]
As usual we are looking for the solutions with $E\in [0,\la]$. Then the functions $y_1$ and $y_2$
are of type \eqref{eqgensolbas} and \eqref{eqgensolcon} (now on intervals $(-R/2,R/2)$, $(-R/2, r+R/2)$).

But it is easy to show that if
\[
A \sin(\sqrt\la (\xi+x))+B \cos(\sqrt\la (\xi+x))
\]
is an even function of $x$, then it is of the form $C \cos(\sqrt\la x)$
(and the same holds for the linear combinations of $\sinh$ and $\cosh$).
Consequently, by the symmetry mentioned above, we can conclude that
\[
y_{ben}=A \cos  (\sqrt{\frac{\la-E}{a}} x), \quad y_{nb}=B \cosh \left(\sqrt{\frac{\mu+E}{b}}(x-(R+r)/2)\right).
\]
Gluing the functions and their derivatives on the interface $x=R/2$ (recall \eqref{eqtickglue}) yields
\[
A \cos  (\sqrt{\frac{\la-E}{a}} \, \frac{R}{2})
=B \cosh \left(\sqrt{\frac{\mu+E}{b}} \, \frac{r}{2}\right),
\]
\[
A a \sqrt{\frac{\la-E}{a}} \sin  (\sqrt{\frac{\la-E}{a}} \, \frac{R}{2})
=B b \sqrt{\frac{\mu+E}{b}} \sinh \left(\sqrt{\frac{\mu+E}{b}} \, \frac{r}{2}\right),
\]
and thus
\[
\sqrt{a(\la-E)} \tan  \left(\sqrt{\frac{\la-E}{a}}\, \frac{R}{2}\right)
=\sqrt{b(\mu+E)} \tanh \left(\sqrt{\frac{\mu+E}{b}} \, \frac{r}{2}\right),
\]
which is the same as \eqref{eqNeu0pr} but with $R/2$ and $r/2$ instead of $R$ ad $r$.
Hence the proof is completed as in Theorem \ref{th1a}.

 \vspace{3mm}
{\it Proof of Theorem \ref{th3}}.

Now a solution $y$ to the eigenvalue problem $Ly=Ey$ is obtained by gluing solutions $y_{odd}$ and $y_{even}$
defined by  formulas like \eqref{eqgensolbas} and \eqref{eqgensolcon} on odd and even subintervals. The key point is, that
under the periodic condition the operator $L$ commutes with the shift operator $Tf(x)=f(x+R+r)$, where the addition
is understood modulo $R+r$. This means that for any eigenfunction $y$ of $L$, the function $Ty$ is also an eigenfunction
with the same eigenvalue. Using the nondegeneracy of eigenvalues of $L$ (like in the proof of the previous theorem)
we can conclude that
$Ty=\om y$ for some $\om$. But $T^k$ is the identity operator, so that $\om^k=1$. Since our eigenfunctions are real
it follows that $\om =\pm 1$. Again the multiplication by $-1$ would contradict the
continuity and thus only $\om=1$ is allowed. Hence we are directly
in the setting of Theorem \ref{th2} and the proof is complete.


 \vspace{3mm}
{\it Proof of Theorem \ref{thvec1} and \ref{thvec2} }.
Under the condition of the Theorem, there exists an invertible matrix $C$ such that $CMC^{-1}=D$ is diagonal with
the elements $\La_j$ on the diagonal. Then in the variables $z=Cy$ all equations and boundary conditions remain the same,
 but with $M$ substituted by $D$. Hence our system decomposes into $n$ independent equations for the coordinates $z^j$ of the
 vector-valued function $z$. For all $z^j$ with $j>1$ there can be no positive eigenvalues, because all operators involved
 are strictly negative. And the equations for $z^1$ (and all boundary conditions) become identical to the equations
 for $y$ of the one-dimensional case, with $\La_1$ instead of $\la$ and the new $\mu$ being equal to $\la+\mu$ of
 the one-dimensional case.

 \vspace{3mm}
{\it Proof of Theorem \ref{thvec3} }. The eigenvalue problem to equation \eqref{eqtickdiffvec1n}
with the Dirichlet boundary condition writes down as
$Ay''+My=Ey$ subject to $y(0)=y(R)=0$. Equivalently this equation rewrites as
\begin{equation}
 \label{eq1thvec3pr}
y''+(A^{-1}M-A^{-1}E)y=0.
\end{equation}
 It is easy to see that the highest real eigenvalue of the matrix $A^{-1}M-A^{-1}E$ decreases with the increase of $E$.
 Hence, according to Theorem \ref{thvec1}, equation \eqref{eq1thvec3pr} is solvable for some positive $E$ if and only if
 the highest real eigenvalue $\La_1$ of the matrix $A^{-1}M$ satisfies  the condition $\La_1> \pi^2/R^2$
 implying the claim of the theorem.

 \vspace{3mm}
{\it Proof of Theorem \ref{thvec4} }.

 Change the variable in \eqref{eqtwopart1} to $z$ so that
to $z_{ben}=A_{ben}y$, $z_{nb}=A_{nb}y$. Then we get

 \begin{equation}
\label{eqtwopart2}
 \begin{aligned}
 & z''_{ben} +(M_{ben}-E)A_{ben}^{-1} z_{ben}=0 \\
 & z''_{nb} +(M_{nb}-E) A_{nb}^{-1} z_{nb}=0,
\end{aligned}
\end{equation}
and the gluing condition for $z$ is the continuity of $z$ and $z'$ on the interfaces.

Let $L(E)$ denote the operator on the r.h.s. of \eqref{eqtwopart2}, so that system  \eqref{eqtwopart2}
can be concisely written as $L(E)z=0$. We are interested in a criterion that ensures that this equation has no solution
with $E\ge 0$. This would be the case, if one could show that $(z,L(E)z)<0$ for all $z$.
But $(z,L(E)z)=(z,L^S(E)z)$, where $L^S(E)=(L(E)+L^T(E))/2$ is the symmetrization. Thus it is sufficient to show that
$L^S(E)$ is a negative symmetric operator for any positive $E$. But $L^S(E)$ decreases with the increase of $E$. Thus
it is sufficient to show that $L=L(0)$ is negative, that is $Lz=Pz$ has no solutions with non-negative $P$.

And this means that there is no solutions to the system
\begin{equation}
\label{eqtwopart3}
 \begin{aligned}
 & z''_{ben} +(M^S_{ben}-P) z_{ben}=0 \\
 & z''_{nb} + (M^S_{nb}-P) z_{nb}=0,
\end{aligned}
\end{equation}
with $P\ge 0$ and $M^S_{ben}=(M_{ben}A_{ben}^{-1}+A_{ben}^{-1}M_{ben}^T)/2$,
$M^S_{nb}=(M_{nb}A_{nb}^{-1}+A_{nb}^{-1}M_{nb}^T)/2$.

It is clear that the required sufficient condition is sufficient to show for the case when all eigenvalues
of $M^S_{ben}$ and $M^S_{nb}$ are different, as the general case is obtained by a straightforward limiting procedure.

For simplicity, let us assume that $M^S_{ben}$ has one positive eigenvalue (that is $k=1$) and all other are negative:
$\la_1>0\ge  \la_2> \cdots > \la_n$.

Clearly $P\ge \la_1$ are impossible in \eqref{eqtwopart3}, so that only $P\in [0, \la_1]$ have to be analyzed.

 Solutions to \eqref{eqtwopart3} are given by the formulas
\[
z_{ben}=A_1 \cos(\sqrt{\la_1-P} x)v_1+\sum_{j=2}^n A_j \cosh(\sqrt{P-\la_j} x)v_j,
\]
\[
z_{nb}=\sum_{j=1}^n B_j \cosh(\sqrt{P-\mu_j} (x-(R+r)/2))w_j,
\]
where
$v_j$ and $w_j$ are orthonormal spectral bases for $M^S_{ben}$ and $M^S_{nb}$.

Let $w_j=\sum_k c_{jk} v_k$ with an orthogonal matrix $C=(c_{ij})$. Then
\[
 z_{nb}=\sum B_k \cosh(\sqrt{P-\mu_k} (x-(R+r)/2))c_{kj}v_j.
\]

Thus gluing $z_{ben}$ and $z_{nb}$ at $x=R/2$ yields
\[
A_1 \cos(\sqrt{\la_1-P} R/2)=\sum_j B_j \cosh(\sqrt{P-\mu_j}r/2)c_{j1},
\]
and
\[
A_i \cosh(\sqrt{P-\la_j} R/2)=\sum_j B_j \cosh(\sqrt{P-\mu_j}r/2)c_{ji}, \quad i>1.
\]

Next, since
\[
\int_{-R/2}^{r+R/2}(z(x),z''(x)) \, dx=-\int_{-R/2}^{r+R/2}(z'(x),z'(x))\, dx,
\]
which follows from the integration by parts and the gluing (and periodic) conditions, we have that,
for a solution to \eqref{eqtwopart3},
\[
0 \le \int_{ben}(z_{ben},(M^S_{ben}-P)z_{ben})\, dx+\int_{nb}(z_{nb}, (M^S_{nb}-P)z_{nb})\, dx
\]
\[
=A_1^2 (\la_1-P)\int_{ben} \cos^2(\sqrt{\la_1-P} x)\, dx
+\sum_{j>1} A_j^2 (\la_j-P)\int_{ben} \cosh^2(\sqrt{P-\la_j} x)\, dx
\]
\[
+\sum_{j \ge 1} B_j^2 (\mu_j-P)\int_{nb} \cosh^2(\sqrt{P-\mu_j} (x-(R+r)/2))\, dx.
\]

Since $\cos (2x)=2\cos^2 x-1$ and $\cosh (2x)=2\cosh^2 x-1$,
\[
\int_{\xi}^{\xi} \cos^2(\al x) \, dx= \xi+\frac12\int_{\xi}^{\xi} \cos(2\al x) \, dx=\xi+\frac{\sin(2\al \xi)}{2\al},
\]
and thus
\begin{equation}
\label{eq}
A_1^2 (\la_1-P)\left(\frac{R}{2}+\frac{\sin(R \sqrt{\la_1-P})}{2\sqrt{\la_1-P}}\right)
\ge   \sum_{j \ge 1} B_j^2 (P-\mu_j)\left(\frac{r}{2}+\frac{\sinh(r \sqrt{P-\mu_j})}{2\sqrt{P-\mu_j}}\right).
\end{equation}

In particular,
\[
A_1^2 \ge B_j^2\frac{P-\mu_j}{\la_1-P}
\frac{\left(\frac{r}{2}+\frac{\sinh(r \sqrt{P-\mu_j})}{2\sqrt{P-\mu_j}}\right)}
{\left(\frac{R}{2}+\frac{\sin(R \sqrt{\la_1-P})}{2\sqrt{\la_1-P}}\right)}
\]
for any $j$.

One could work with this estimate, but it makes the final result more transparent (though a bit more rough), if
one simplifies  this estimate to the following one:
\[
A_1^2 \ge B_j^2 \frac{\sqrt{P-\mu_j}\sinh(r \sqrt{P-\mu_j})}
{2R(\la_1-P)}.
\]

On the other hand,
\[
|A_1| \cos(\sqrt{\la_1-P} \frac{R}{2}) \le \sum_j |B_j| \cosh(\sqrt{P-\mu_j}\frac{r}{2})c_{j1},
\]
and thus
\[
A_1^2 \cos^2(\sqrt{\la_1-P} \frac{R}{2}) \le n\sum_j B_j^2 \cosh^2(\sqrt{P-\mu_j}\frac{r}{2})c^2_{j1}
\le n \max_j B_j^2 \cosh^2(\sqrt{P-\mu_j}\frac{r}{2}),
\]
or
\[
A_1^2 (1+\cos(\sqrt{\la_1-P}R))
\le n \max_j B_j^2(1+ \cosh(\sqrt{P-\mu_j}r)).
\]

Thus for a $j$ realising the maximum we have

\[
n \frac{1+\cosh(\sqrt{P-\mu_j}r)}{1+\cos(\sqrt{\la_1-P} R)}
 \ge
\frac{\sqrt{P-\mu_j}\sinh(r \sqrt{P-\mu_j})}
{2R(\la_1-P)},
\]
or
 \begin{equation}
\label{eqcondviasym4}
\frac{R(\la_1-P)}{(1+\cos(R\sqrt{\la_1-P}))}
 \ge
\frac{\sqrt{P-\mu_j}\sinh(r \sqrt{P-\mu_j})}
{2n(1+\cosh(r \sqrt{P-\mu_j}))}.
\end{equation}
This can be rewritten as
 \begin{equation}
\label{eqcondviasym5}
\frac{1}{R} g_1(R\sqrt{\la_1-P})
 \ge \frac{1}{2nr} g_2(r\sqrt{P-\mu_j})
\end{equation}
with
 \begin{equation}
\label{eqcondviasym6}
g_1(x)= \frac{x^2}{1+\cos x}, \quad
g_2(x)= \frac{x\sinh x}{1+\cosh x}.
\end{equation}

In order to make \eqref{eqcondviasym5} impossible (and thus to ensure the negativity of our operator)
it is sufficient to assume that
 \begin{equation}
\label{eqcondviasym7}
\frac{1}{R} g_1(R\sqrt{x})
 < \frac{1}{2nr} g_2(r\sqrt{\la_1-\mu_1-x})
\end{equation}
for $x\in [0, \la_1]$, where $x=\la_1-P$.   For $\la_1 \ge \pi^2/R$, this cannot be true. Hence assume
$\la_1 < \pi^2/R$.

Since it is straightforward to check that the function
 \[
\frac{ g_1(R\sqrt{x})}{g_2(r\sqrt{\la_1-\mu_1-x})}
\]
is increasing for $x\in [0, \la_1]$, condition \eqref{eqcondviasym7} is equivalent
to this condition on the endpoint:
   \begin{equation}
\label{eqcondviasym8}
\frac{1}{R} g_1(R\sqrt{\la_1})
 < \frac{1}{2nr} g_2(r\sqrt{|\mu_1|}),
\end{equation}
or explicitly
    \begin{equation}
\label{eqcondviasym9}
\frac{2Rn \la_1}{1+\cos (R\sqrt{\la_1}))}
 < \frac{\sqrt{|\mu_1|} \sinh (r\sqrt{|\mu_1|})}
 {1+\cosh(r\sqrt{|\mu_1|})},
\end{equation}
which coincides with \eqref{eq1thvec4} for $k=1$.

Modifications required for arbitrary $k$ are straightforward, and we omit them.

 \vspace{3mm}
{\it Proof of Theorem \ref{thvec5} }.

The eigenvalue equations \eqref{eq2dimvarst} can be written as
 \begin{equation}
\label{eq2dimvarst1}
 \begin{aligned}
 & y''_{ben} +N_{ben}(E)y_{ben}=0 \\
 & y''_{nb} +N_{nb}(E)y_{nb}=0.
\end{aligned}
\end{equation}

As above we can write
\[
y_{ben}=A\cos (\sqrt{\La_1} x)v_1+B\cosh (\sqrt{|\La_2|} x)v_2,
\]
\[
y_{nb}=C\cosh (\sqrt{|\mu_1|} (x-(R+r)/2))w_1+D\cosh (\sqrt{|\mu_2|} (x-(R+r)/2))w_2.
\]
The last equation rewrites as
\[
y_{nb}=C\cosh (\sqrt{|\mu_1|} (x-(R+r)/2))(c_{11}v_1+c_{12}v_2)+D\cosh (\sqrt{|\mu_2|} (x-(R+r)/2))(c_{21}v_1+c_{22}v_2).
\]
Of course, all terms here depend on $E$.

Further on, in order to shorten the formulas, we shall write $\sin$, $\cos $, $\tan$ if the argument is $\sqrt{\La_1} R/2$
 and $\sinh$, $\cosh $, $\tanh$ if the argument is $\sqrt{|\La_2|} R/2$.
 Gluing the solutions at the interface $x=R$ yields
\[
A\cos = C\cosh (\sqrt{|\mu_1|}r/2) c_{11}+D\cosh (\sqrt{|\mu_2|} r/2)c_{21},
\]
\[
B\cosh =C\cosh (\sqrt{|\mu_1|} r/2) c_{12}+D\cosh (\sqrt{|\mu_2|} r/2)c_{22}.
\]
and
\[
-A\sqrt{\La_1} \sin
=-a(C\sqrt{|\mu_1|} \sinh (\sqrt{|\mu_1|} r/2)c_{11}+D\sqrt{|\mu_2|}\sinh (\sqrt{|\mu_2|} r/2)c_{21}),
\]
\[
B\sqrt{\La_1} \sinh
=-a(C\sqrt{|\mu_1|} \sinh (\sqrt{|\mu_1|} r/2)c_{12}+D\sqrt{|\mu_2|}\sinh (\sqrt{|\mu_2|} r/2)c_{22}).
\]

Equalising $A$ from  the 1st equations yields
\[
\sqrt{\La_1} \tan [C\cosh (\sqrt{|\mu_1|}r/2) c_{11}+D\cosh (\sqrt{|\mu_2|} r/2)c_{21}]
\]
\[
= a[C\sqrt{|\mu_1|} \sinh (\sqrt{|\mu_1|} r/2)c_{11}+D\sqrt{|\mu_2|}\sinh (\sqrt{|\mu_2|} r/2)c_{21}],
\]
Equalising $B$ from  the 2nd equations yields
\[
\sqrt{|\La_2|} \tanh [C\cosh (\sqrt{|\mu_1|} r/2) c_{12}+D\cosh (\sqrt{|\mu_2|} r/2)c_{22}]
\]
\[
=- a[ C\sqrt{|\mu_1|} \sinh (\sqrt{|\mu_1|} r/2)c_{12}+D\sqrt{|\mu_2|}\sinh (\sqrt{|\mu_2|} r/2)c_{22}].
\]
The condition for the existence of a solution to this system of two equations with two unknown $C,D$
writes down as  the equation
\[
c_{11}c_{22} [\sqrt {\La_1} \tan -a\sqrt {|\mu_1|} \tanh (\sqrt{|\mu_1|}r/2)]
[\sqrt {|\La_2|} \tanh +a\sqrt {|\mu_2|} \tanh (\sqrt{|\mu_2|}r/2)]
\]
\[
-c_{12}c_{21} [\sqrt {\La_1} \tan -a\sqrt {|\mu_2|} \tanh (\sqrt{|\mu_2|}r/2)]
[\sqrt {|\La_2|} \tanh +a\sqrt {|\mu_1|} \tanh (\sqrt{|\mu_1|}r/2)]=0,
\]
or as $\sqrt \La_1 \tan =aB/A$ with
\[
A=c_{11}c_{22}
[\sqrt {|\La_2|} \tanh +a\sqrt {|\mu_2|} \tanh (\sqrt{|\mu_2|}r/2)]
\]
\[
-c_{12}c_{21}
[\sqrt {|\La_2|} \tanh +a\sqrt {|\mu_1|} \tanh (\sqrt{|\mu_1|}r/2)],
\]
\[
B=c_{11}c_{22}\sqrt {|\mu_1|} \tanh (\sqrt{|\mu_1|}r/2)
[\sqrt {|\La_2|} \tanh +a\sqrt {|\mu_2|} \tanh (\sqrt{|\mu_2|}r/2)]
\]
\[
-c_{12}c_{21} \sqrt {|\mu_2|} \tanh (\sqrt{|\mu_2|}r/2)]
[\sqrt {|\La_2|} \tanh +a\sqrt {|\mu_1|} \tanh (\sqrt{|\mu_1|}r/2)].
\]

Since $c_{12}c_{21} \le 0$ and $c_{11}c_{22} \ge 0$ (and $|\mu_2|\ge |\mu_1|)$, it follows that
\[
B\ge \det c [a\sqrt {|\mu_2|} \tanh (\sqrt{|\mu_2|}r/2)+\sqrt {\La_2} \tanh] \sqrt {|\mu_1|} \tanh (\sqrt{|\mu_1|}r/2),
\]
\[
0<A< \det c [\sqrt {|\La_2|} \tanh +a\sqrt {|\mu_2|} \tanh (\sqrt{|\mu_2|}r/2)],
\]
so that
\begin{equation}
\label{eqestimAB0}
a B/ A\ge a\sqrt {|\mu_1|} \tanh (\sqrt{|\mu_1|}r/2).
\end{equation}
Hence
\[
\sqrt \La_1 \tan (\sqrt{\la_1}R/2) \ge
a\sqrt {|\mu_1|} \tanh (\sqrt{|\mu_1|}r/2).
\]
Consequently, if \eqref{eq1thvec5} holds,
this is impossible. It remains to notice that,
by monotonicity, this condition holds for all $E$ if and only if it holds for $E=0$.

 \vspace{3mm}
{\it Proof of Theorem \ref{thvec6} }.

As is easily seen, if inequalities \eqref{eq1thvec6} hold for $N_{ben}(0)$ and $N_{nb}(0)$, then
they hold for $N_{ben}(E)$ and $N_{nb}(E)$ with all positive $E$. Hence it is sufficient to check
condition \eqref{eq2thvec5} for $E=0$.

For a matrix
\[
A=\left(\begin{aligned}
& -m_1 \quad b_1 \\
& \,\, b_2 \quad -m_2
\end{aligned}
\right),
\]
the eigenvalues are
\[
\La_{1,2}=\frac12\left[ -(m_1+m_2) \pm \sqrt{(m_1-m_2)^2+4b_1b_2}\right].
\]
Thus $\La_2<0$ and $\La_1>0$.

Looking for eigenvectors $v_{1,2}$, let us choose $v_1^1=1, v_2^2=1$. Then
\[
v_1^2=\frac{\La_1+m_1}{b_1}=\frac{b_2}{\La_1+m_2},
\]
\[
v_2^1=  \frac{b_1}{\La_2+m_1}=\frac{\La_2+m_2}{b_2}.
\]

Thus
\[
V=\left(
\begin{aligned}
& 1 \quad v_2^1 \\
& v_1^2 \quad 1
\end{aligned}\right),
\quad
W=\left(
\begin{aligned}
& 1 \quad w_2^1 \\
& w_1^2 \quad 1
\end{aligned}\right).
\]
Consequently,
\[
C=WV^{-1}
=\frac{1}{\det V}\left(
\begin{aligned}
& 1 -v_2^1 w_1^2\quad w_2^1-v_2^1 \\
& -v_1^2 +w_1^2\quad 1-v_1^2w_2^1
\end{aligned}\right).
\]

Assuming $m_1>m_2$, we see that
 \[
 v_2^1=-\frac{(m_1-m_2)+\sqrt{(m_1-m_2)^2+4b_1b_2}}{2b_2} <0,
 \]
 \[
 v_1^2=\frac{(m_1-m_2)+\sqrt{(m_1-m_2)^2+4b_1b_2}}{2b_1} <0.
 \]

Similarly, $w_2^1<0$ and $ w_1^2>0$. They have the same expression, but with
$\tilde m_j$ and $\tilde b_j$ instead of $m_j,b_j$.
Hence $c_{11}, c_{22} \ge 1>0$. Moreover,
\[
c_{12}c_{21}=-(w_1^2-v_1^2)(|w_2^1|-|v_2^1|).
\]
Under \eqref{eq0thvec6} and \eqref{eq00thvec6}, it follows that, for any $\om \in (0,1)$,
\[
w_1^2 \ge v_1^2, |w_2^1| \ge |v_2^1|,
\]
and thus $c_{12}c_{21} <0$, which completes the proof.

\section{Conclusion}

In this paper we have developed the diffusion model for the propagation
of ticks with variable discontinuous diffusion coefficients. We have extended the famous KISS model
(originated in plankton research) to describe the dynamics on patchy territories with intermittent background
(beneficial) and control (non-beneficial) zones and gave conditions for the control parameters that ensure
the eradication of ticks' population. Mathematically these conditions imply the negativity of the corresponding
diffusion operator. We have started with rather trivial (though complete) results for the scalar case
(one-stage populations) and then presented some advanced theorems for the vector-valued non-symmetric
diffusions stressing serious theoretical difficulties that arise in this extension. Finally we have used various
published resources reporting on the concrete observations of ticks in order to extract some concrete numeric values
for the parameters of background zones and thus to provide approximate numeric values to the required parameters
of the control zones. In particular, we found that the most rough scalar model yields closer values of about 16 km for the KISS
critical patch sizes related to both the lone star and the taiga ticks. On the other hand, using more detailed information combined
with a more detailed multi-stage model can provide quite different output parameters. Also concrete numeric data
showed that the estimates based on the symmetrization of the birth-and-death matrix (vary natural from the
theoretical point of view) can underestimate strongly the lower  bounds for the KISS scales.

Let us note some further perspectives.
We used in this paper a one-dimensional diffusion models aiming at the approximate description of the propagation in a
certain direction. Of course, it would be natural to extend the results to a more realistic (for ticks and many other species)
two-dimensional diffusion models. On the other hand, we discussed the solutions in a global closed territory. Of interest would
be a more local analysis of the solutions in a (possibly intermittent) small territory placed in contact with a large reservoir
of an expanding background population. Such solutions could be possibly analysed based on the Fokas method, see \cite{Fokas}
and \cite{SheDec}, or via the method of estimating the heat content, see e.g. \cite{VanBerg15}.
It would be also of interest to see whether the results of Theorem \ref{thvec5} have any reasonable extensions
beyond the restrictive assumptions \eqref{eq2thvec5} and whether one can get some multi-dimensional analogs.
Also Remark \ref{remimprovesym} can be used for practical calculations with multi-dimensional models. A realistic
vector-valued diffusion model can include not only ticks in various stages, but also the hosts.

Practically, the main goal of this paper was to give exact recommendations on the properties of control
zones needed for the eradication of ticks. But the numbers obtained in our analysis are based on very rough estimates,
and our calculations show that they are quite sensitive to the input data.
Much more detailed program of concrete observations in various places is needed to make the estimates more precise
and more specific for concrete regions.

In this paper we stick to the diffusion models. Further
investigations are needed to see, whether more appropriate approximations can be achieved by using other models,
like the telegraph equations, or presently very popular sub-diffusions or super-diffusions.

Notice finally that our models with periodic intermittent beneficial/non-beneficial zones present much
simplified versions of the parabolic Anderson models with random potentials, to which an immense amount
of studies were devoted, see e.g. \cite{MHMY} and references therein.

{\bf Acknowledgements.} This work was supported by
 the Ministry of Science and Higher Education of the Russian Federation Grant ID: 075-15-2020-928.

\end{document}